# QUASIHYPERBOLIC GEODESICS IN JOHN DOMAINS IN $\mathbb{R}^n$

M. HUANG AND X. WANG *

ABSTRACT. In this paper, we prove that if $D \subset \mathbb{R}^n$ is a John domain which is homeomorphic to a uniform domain via a quasiconformal mapping, then each quasihyperbolic geodesic in $D$ is a cone arc, which shows that the answer to one of open problems raised by Heinonen in [6] is affirmative. This result also shows that the answer to the open problem raised by Gehring, Hag and Martio in [5] is positive for John domains which are homeomorphic to uniform domains via quasiconformal mappings. As an application, we prove that if $D \subset \mathbb{R}^n$ is a John domain which is homeomorphic to a uniform domain, then $D$ must be a quasihyperbolic $(b, \lambda)$-uniform domain.

## 1. INTRODUCTION AND THE MAIN RESULT

In the following, we always assume that $D$ is a proper subdomain in $\mathbb{R}^n$. We begin with the following concepts.

**Definition 1.1.** $D$ is called *c-uniform* if there exists a constant $c$ with the property that each pair of points $z_1, z_2$ in $D$ can be joined by a rectifiable arc $\gamma$ in $D$ satisfying (cf. [10, 14])

(1) $\min_{j=1,2} \ell(\gamma[z_j, z]) \le c\, d(z)$ for all $z \in \gamma$, and

(2) $\ell(\gamma) \le c\, |z_1 - z_2|$,

where $\ell(\gamma)$ denotes the arclength of $\gamma$, $\gamma[z_j, z]$ the part of $\gamma$ between $z_j$ and $z$, and $d(z)$ the distance from $z$ to the boundary $\partial D$ of $D$. Also we say that $\gamma$ is a *double c-cone arc.*

$D$ is said to be a *c-John domain* if it satisfies the condition (1) in Definition 1.1, not necessarily (2), and $\gamma$ is called a *c-cone arc.*

John [12], Martio and Sarvas [10] were the first who introduced John domains and uniform domains, respectively. Now, there are plenty of alternative characterizations for uniform and John domains (see [1, 2, 3, 8, 9]). And its importance along with some special domains throughout the function theory is well documented, see [2, 8, 11].

Gehring and Osgood [3] proved that each quasihyperbolic geodesic in a $c$-uniform domain $D \subset \mathbb{R}^n$ is a double $b$-cone arc, where the constant $b$ depends only on $c$. Since a John domain can be thought as a “one-sided” uniform domain, it is natural

---

* Corresponding author.

The research was partly supported by NSFs of China (No. 10771059 and 11071063).

to ask whether the result is true or not for John domains. In fact, this problem has been proposed by Gehring, Hag and Martio in [5] in the following form.

**Conjecture 1.2.** Suppose that $D \subset \mathbb{R}^n$ is a $c$-John domain with center $x_0$ and that $\gamma$ is a quasihyperbolic geodesic which joins $x_1$ to $x_0$. Is $\gamma$ a $b$-cone arc for some $b = b(c)$?

Gehring, Hag and Martio themselves discussed Conjecture 1.2 and got the following.

**Theorem A.** ([5, Theorem 4.1]) *Suppose that $D \subset \mathbb{R}^2$ is a simply connected $c$-John domain. If $\gamma$ is either a quasihyperbolic or hyperbolic geodesic in $D$, then $\gamma$ is a $b$-cone arc, where $b$ depends only on $c$.*

Theorem A shows that the answer to Conjecture 1.2 is yes when $n = 2$ and $D$ is simply connected. In [5], the authors also constructed counterexamples to show that the answer to Conjecture 1.2 is no when $D$ is multiply connected or $D$ is simply connected and $n > 2$. These counterexamples explain that to study Conjecture 1.2 further, some restriction is necessary. Hence, in [6], Heinonen modified Conjecture 1.2 to the following form.

**Conjecture 1.3.** Suppose that $D \subset \mathbb{R}^n$ is a $c$-John domain with center $x_0$ and that $D$ is quasiconformally equivalent to the unit ball $\mathbb{B}^n \subset \mathbb{R}^n$. Let $\gamma$ be a quasihyperbolic geodesic which joins $x_1$ to $x_0$. Is $\gamma$ a $b$-cone arc for some $b = b(c)$?

The main aim of this paper is to discuss Conjectures 1.2 and 1.3. We get the following result whose proof will be presented in Section 3.

**Theorem 1.4.** *Suppose that $D \subset \mathbb{R}^n$ is an $a$-John domain which is homeomorphic to a $c$-uniform domain via a $K$-quasiconformal mapping $f$. Let $z_1$, $z_2 \in D$ and $\gamma$ be a quasihyperbolic geodesic joining $z_1$ and $z_2$ in $D$. Then $\gamma$ is an $a'$-cone arc, where the positive constant $a'$ depends only on $a$, $c$, $n$ and $K$.*

**Remark 1.5.** (*i*) *Theorem 1.4 shows that the answer to Conjecture 1.3 is positive. In fact, we have proved more than is stated in Conjecture 1.3 because* (1) *our result is independent of the center $x_0$;* (2) *the condition "the unit ball" in Conjecture 1.3 is replaced by the one "uniform domains". It is known that the unit ball is a $\frac{\pi}{2}$-unform domain.*

(*ii*) (1) *Theorem 1.4 also shows that the answer to Conjecture 1.2 is positive for John domains which are homeomorphic to uniform domains via quasiconformal mappings.* (2) *Even when $n = 2$, Riemann mapping theorem shows Theorem 1.4 is a generalization of Theorem A.*

In [8], Kim and Langmeyer got the following result concerning the quasihyperbolic $(b, \lambda)$-uniform domains (see Section 2 for the definition).

**Theorem B.** ([8, Theorem 4.1 and Remark 4.14]) *Suppose $D \subset \mathbb{R}^n$ is an $a$-John domain which is a $K$-quasiconformal image of a $c$-uniform domain in $\mathbb{R}^n$. If each quasihyperbolic geodesic in $D$ is an inner $c'$-cone arc, where $c'$ depends on $a$ and $c$,* *then $D$ is a quasihyperbolic $(b, \lambda)$-uniform domain. Here $b$ depends on the constant $c'$, $K$ and $n$.*

By using Theorem 1.4, in Section 4, we will prove the following result.

**Theorem 1.6.** *Suppose that $D \subset \mathbb{R}^n$ is an $a$-John domain which is homeomorphic to a $c$-uniform domain via a $K$-quasiconformal $f$. Then $D$ is a quasihyperbolic $(b, \lambda)$-uniform domain, where the positive constant $b$ depends only on $a$, $c$, $K$ and $n$.*

**Remark 1.7.** *Theorem 1.6 is a generalization of Theorem B since Theorem 1.6 shows that the hypothesis "each quasihyperbolic geodesic in $D$ being an inner $c'$-cone arc" in Theorem B is redundant.*

## 2. Preliminaries

Let $\gamma$ be a rectifiable arc or a path in $D$. Then the *quasihyperbolic length* of $\gamma$ is the number (cf. [4]):

$$\ell_k(\gamma) = \int_\gamma \frac{|dz|}{d(z)}.$$

For any $z_1$, $z_2$ in $D$, the *quasihyperbolic distance* $k_D(z_1, z_2)$ between $z_1$ and $z_2$ is defined in the usual way:

$$k_D(z_1, z_2) = \inf \ell_k(\gamma),$$

where the infimum is taken over all rectifiable arcs $\gamma$ joining $z_1$ to $z_2$ in $D$. For any $z_1$, $z_2$ in $D$, we have (cf. [3])

$$k_D(z_1, z_2) \geq \log\Big(1 + \frac{|z_1 - z_2|}{\min\{d(z_1), d(z_2)\}}\Big) \geq \Big|\log \frac{d(z_2)}{d(z_1)}\Big|.$$

As a generalization of quasiconformal mappings, Väisälä introduced CQH homeomorphisms (cf. [14]).

**Definition 2.1.** Suppose $f : D \mapsto D'$ is a homeomorphism. Then $f$ is said to be *$C$-coarsely $M$-quasihyperbolic*, or briefly $(M, C)$-CQH, in the quasihyperbolic metric if it satisfies

$$\frac{k_D(x, y) - C}{M} \leq k_{D'}(f(x), f(y)) \leq M\ k_D(x, y) + C$$

for all $x$, $y \in D$.

The following proposition easily follows from [3, Theorem 3].

**Proposition 2.2.** *Each $K$-quasiconformal mapping in $\mathbb{R}^n$ is an $(M, C)$-CQH homeomorphism with $(M, C)$ depending only on $(K, n)$.*

See [13, 16] for more details about the properties of quasiconformal mappings.

Let's recall the following characterization of uniform domains, which is due to Gehring and Osgood.

**Theorem C.** ([3, Theorem 1]) *$D$ is a $c$-uniform domain if and only if $k_D(x,y) \leq c' \log\left(1+\frac{|x-y|}{\min\{d(x),d(y)\}}\right)$ for any $x, y \in D$, where the constants $c$ and $c'$ depend only on each other.*

For any $z_1$, $z_2 \in D$, the *inner distance* $\lambda_D(z_1, z_2)$ between them is defined by

$$\lambda_D(z_1, z_2) = \inf\{\ell(\alpha) : \ \alpha \subset D \text{ is a rectifiable arc joining } z_1 \text{ and } z_2\}.$$

**Definition 2.3.** $D$ is called *inner $c$-uniform* if there exists a constant $c$ with the property that each pair of points $z_1, z_2$ in $D$ can be joined by a rectifiable arc $\gamma$ in $D$ satisfying (cf. [15])

(1) $\min\limits_{j=1,2} \ell(\gamma[z_j, z]) \leq c\, d(z)$ for all $z \in \gamma$, and

(2) $\ell(\gamma) \leq c\, \lambda_D(z_1, z_2)$.

Also we say that $\gamma$ is an *inner double $c$-cone arc.*

Väisälä introduced the concept of quasihyperbolic $(b, \lambda)$-uniform domain in [15].

**Definition 2.4.** A domain $D$ in $\mathbb{R}^n$ is said to be a *quasihyperbolic $(b, \lambda)$-uniform domain*, or briefly *QH $(b, \lambda)$-uniform* if

$$k_D(z_1, z_2) \leq b \log\left(1 + \frac{\lambda_D(z_1, z_2)}{\min\{d(z_1), d(z_2)\}}\right)$$

for all $z_1, z_2 \in D$, where $b \geq 1$ is a constant.

Obviously, inner $c$-uniform domains and $QH$ $(b, \lambda)$-uniform domains are generalizations of uniform domains. The following result describes the relation between inner uniform domains and $QH$ $(b, \lambda)$-uniform domains.

**Theorem D.** ([15, Theorem 2.33]) *$D$ is an inner $c$-uniform domain if and only if $D$ is a QH $(b, \lambda)$-uniform domain, where the constants $b$ and $c$ depend on each other.*

## 3. The proof of Theorem 1.4

In what follows, we always assume that $f : D \mapsto D'$ is a $K$-quasiconformal mapping. Also we use $x$, $y$, $z$, $\cdots$ to denote the points in $D$, and $x'$, $y'$, $z'$, $\cdots$ the images of $x$, $y$, $z$, $\cdots$ in $D'$, respectively, under $f$. For arcs $\alpha$, $\beta$, $\gamma$, $\cdots$ in $D$, we also use $\alpha'$, $\beta'$, $\gamma'$, $\cdots$ to denote their images in $D'$.

For $x$, $y \in D$ , let $\beta$ be an arc joining $x$ and $y$ in $D$. We come to determine some special points on $\beta'$.

3.1. **Determination of special points on $\beta'$.** Without loss of generality, we may assume that $d(y') \geq d(x')$. Then there must exist a point $w_0' \in \beta'$ which is the first point along the direction from $x'$ to $y'$ such that

$$d(w_0') = \sup_{p \in \beta'} d(p).$$

It is possible that $w_0' = x'$ or $y'$. Obviously, there exists a nonnegative integer $m$ such that

$$2^m\, d(x') \le d(w_0') < 2^{m+1}\, d(x'),$$

and $x_0'$ the first point in $\beta'[x', w_0']$ from $x'$ to $w_0'$ with

$$d(x_0') = 2^m\, d(x').$$

Let $x_1' = x'$. If $x_0' = x_1'$, we let $x_2' = w_0'$. It is possible that $x_1' = x_2'$. If $x_0' \ne x_1'$, then we let $x_2', \dots, x_{m+1}' \in \beta'[x', x_0']$ be the points such that for each $i \in \{2, \dots, m+1\}$, $x_i'$ denotes the first point from $x'$ to $x_0'$ with

$$d(x_i') = 2^{i-1}\, d(x_1').$$

Obviously, $x_{m+1}' = x_0'$. If $x_0' \ne w_0'$, then we use $x_{m+2}'$ to denote $w_0'$.

In a similar way, let $s \ge 0$ be the integer such that

$$2^s\, d(y') \le d(w_0') < 2^{s+1}\, d(y'),$$

and $x_{1,0}'$ the first point in $\beta'[y', x_{1,0}']$ from $y'$ to $x_{1,0}'$ with

$$d(x_{1,0}') = 2^s\, d(y').$$

Let $x_{1,1}' = y'$. If $x_{1,0}' = x_{1,1}'$, we let $x_{1,2}' = x_{1,0}'$. It is possible that $x_{1,2}' = x_{1,1}'$. If $x_{1,0}' \ne y'$, then we let $x_{1,2}', \dots, x_{1,s+1}'$ be the points in $\beta'[y', w_0']$ such that for each $j \in \{2, \dots, s+1\}$, $x_{1,j}'$ is the first point from $x_{1,1}'$ to $w_0'$ with

$$d(x_{1,j}') = 2^{j-1}\, d(x_{1,1}').$$

Then $x_{1,s+1}' = x_{1,0}'$. If $x_{1,0}' \ne w_0'$, we let $x_{1,s+2}' = w_0'$.

3.2. **Elementary properties.** In the following, we assume that for any $s_1$, $s_2 \in \beta$,

$$\ell_k(\beta[s_1, s_2]) \le 4a^2 k_D(s_1, s_2) + 4a^2. \tag{3.1}$$

Obviously, for each quasihyperbolic geodesic, (3.1) is satisfied.

By Proposition 2.2, in the following, we assume that $f : D \mapsto D'$ is an $(M, C)$-CQH homeomorphism, where $(M, C)$ depends only on $(K, n)$.

**Lemma 3.2.** *For any $k \in \{1, \cdots, m\}$ and $z' \in \beta'[x_k', x_{k+1}']$,*

(1) $d(x_{k+1}') \le a_2\, d(z')$;

(2) $|x_{k+1}' - x_k'| \le a_2\, d(z')$ *and*

(3) $\max\{|x_k' - z'|, |x_{k+1}' - z'|\} \le a_2\, d(z')$,

*where $a_2 = (1 + 2a_1)^{4a^2c'M^2+1} e^{C+4a^2M+4a^2CM}$, $a_1 = e^{3(C+1)(a_0+M)}$ and $a_0 = 2^4[c' + 4a^2c'M + C + 4a^2]^4$. Here and in what follows, $[\cdot]$ always denotes the greatest integer part.*

**Proof.** At first, we prove the following inequality: For any $k \in \{1, \cdots, m\}$,

$$|x_{k+1}' - x_k'| < a_1\, d(x_{k+1}'). \tag{3.3}$$

We prove this inequality by contradiction. Suppose

$$|x_{k+1}' - x_k'| \ge a_1\, d(x_{k+1}'). \tag{3.4}$$

Let $y'_{k,1}, y'_{k,2}, \cdots, y'_{k,a_0+1} \in \beta'[x'_k, x'_{k+1}]$ be $a_0+1$ points such that $y'_{k,1} = x'_k$, $y'_{k,a_0+1} = x'_{k+1}$ and $|y'_{k,i+1} - y'_{k,i}| \geq \frac{|x'_k - x'_{k+1}|}{a_0}$. Then for each $i \in \{1, 2, \cdots, a_0\}$,

$$\begin{aligned} k_{D'}(y'_{k,i}, y'_{k,i+1}) &\geq \log\Big(1 + \frac{|y'_{k,i+1} - y'_{k,i}|}{\min\{d(y'_{k,i+1}), d(y'_{k,i})\}}\Big) \\ &\geq \log\Big(1 + \frac{|x'_k - x'_{k+1}|}{2a_0 d(x'_k)}\Big). \end{aligned}$$

We see from (3.1) and Theorem C that

$$\begin{aligned} a_0 \log\Big(1 + \frac{|x'_k - x'_{k+1}|}{2a_0 d(x'_k)}\Big) &\leq \sum_{i=1}^{a_0} k_{D'}(y'_{k,i}, y'_{k,i+1}) \\ &\leq M \sum_{i=1}^{a_0} k_D(y_{k,i}, y_{k,i+1}) + a_0 C \\ &\leq M\ell_k(\beta[x_k, x_{k+1}]) + a_0 C \\ &\leq 4a^2 M k_D(x_k, x_{k+1}) + 4a^2 M + a_0 C \\ &\leq 4a^2 M^2 k_{D'}(x'_k, x'_{k+1}) + (a_0 + 4a^2 M)C + 4a^2 M \\ &\leq 4a^2 c' M^2 \log\Big(1 + \frac{|x'_k - x'_{k+1}|}{d(x'_k)}\Big) + (a_0 + 4a^2 M)C + 4a^2 M, \end{aligned}$$

whence

$$a_0 \log\Big(1 + \frac{|x'_{k+1} - x'_k|}{2a_0 d(x'_k)}\Big) \leq 8a^2 c' M^2 \log\Big(1 + \frac{|x'_{k+1} - x'_k|}{d(x'_k)}\Big),$$

which contradicts with (3.4). Hence (3.3) holds.

We infer from (3.3) that for any $z' \in \beta'[x'_k, x'_{k+1}]$,

$$\begin{aligned} \text{(3.5)} \quad \log \frac{d(x'_{k+1})}{d(z')} &< k_{D'}(z', x'_{k+1}) \\ &\leq M k_D(z, x_{k+1}) + C \\ &\leq 4a^2 M \; k_D(x_k, x_{k+1}) + C + 4a^2 M \\ &\leq 4a^2 M^2 \; k_{D'}(x'_k, x'_{k+1}) + 4a^2 CM + C + 4a^2 M \\ &\leq 4a^2 c' M^2 \; \log\Big(1 + \frac{|x'_{k+1} - x'_k|}{d(x'_k)}\Big) + C + 4a^2 CM + 4a^2 M \\ &\leq 4a^2 c' M^2 \; \log(1 + 2a_1) + C + 4a^2 CM + 4a^2 M, \end{aligned}$$

which implies that Lemma 3.2 (1) holds. (3.3) and (3.5) yield that

$$|x'_{k+1} - x'_k| \leq (1 + 2a_1)^{4a^2 c' M^2 + 1} e^{4a^2 CM + 4a^2 M + C} \; d(z'),$$

whence Lemma 3.2 (2) follows.

Obviously,

$$\begin{aligned}\log\Big(1+\frac{|x'_k-z'|}{d(z')}\Big) &\leq k_{D'}(x'_k,z')\\ &\leq Mk_D(x_k,z)+C\\ &\leq 4a^2Mk_D(x_k,x_{k+1})+4a^2M+C\\ &\leq 4a^2M^2\ k_{D'}(x'_k,x'_{k+1})+C+4a^2M+4a^2CM\\ &\leq 4a^2c'M^2\ \log\Big(1+\frac{|x'_{k+1}-x'_k|}{d(x'_k)}\Big)+C+4a^2M+4a^2CM,\end{aligned}$$

which, together with (3.3), yields

$$|x'_k-z'|\leq(1+2a_1)^{4a^2c'M^2}e^{C+4a^2M+4a^2CM}\ d(z'). \tag{3.6}$$

The similar discussion as in (3.6) shows that

$$|x'_{k+1}-z'|\leq(1+2a_1)^{4a^2c'M^2}e^{C+4a^2M+4a^2CM}\ d(z'). \tag{3.7}$$

The combination of (3.6) and (3.7) shows that Lemma 3.2 (3) holds. □

The following two results easily follow from the similar reasoning as in the proof of Lemma 3.2.

**Corollary 3.8.** *For any $k\in\{1,\cdots,s\}$ and $z'\in\beta'[x'_{1,k},x'_{1,k+1}]$,*

(1) $d(x'_{1,k+1})\leq a_2\ d(z')$;
(2) $|x'_{1,k+1}-x'_{1,k}|\leq a_2\ d(z')$ *and*
(3) $\max\{|x'_{1,k}-z'|,|x'_{1,k+1}-z'|\}\leq a_2\ d(z')$.

**Corollary 3.9.** *For any $z'\in\beta'[x'_{m+1},x'_{1,s+1}]$,*

(1) $d(w'_0)\leq a_2\ d(z')$;
(2) $|x'_{m+1}-x'_{1,s+1}|\leq a_2\ d(z')$ *and*
(3) $\max\{|x'_{m+1}-z'|,|x'_{1,s+1}-z'|\}\leq a_2\ d(z')$.

**Lemma 3.10.** *For any $z'\in\beta'[x',w'_0]$,*

$$|x'-z'|\leq a_3\ d(z'),$$

*where $a_3=a_2+a_2^2$.*

**Proof.** If $z'\in\beta'[x',x'_{m+1}]$, then there exists some $k\in\{1,\cdots,m\}$ such that $z'\in\beta'[x'_k,x'_{k+1}]$. If $k=1$, then the result easily follows from Lemma 3.2. If $k>1$, then by Lemma 3.2,

$$\begin{aligned}|x'-z'| &\leq |x'_1-x'_2|+\cdots+|x'_{k-1}-x'_k|+|x'_k-z'|\\ &\leq a_2\big(d(x'_1)+\cdots+d(x'_{k-1})+d(z')\big)\\ &\leq (a_2+\frac{1}{2}a_2^2)d(z').\end{aligned}$$

Now we consider the case $z' \in \beta'[x'_{m+1}, w'_0]$. Then we infer from Lemma 3.2 and Corollary 3.9 that

$$\begin{aligned} |x'-z'| &\leq a_2\big(d(x'_1)+d(x'_2)+\cdots+d(x'_m)+d(z')\big)\\ &\leq a_2\big(d(x'_{m+1})+d(z')\big)\\ &\leq (a_2+a_2^2)d(z'). \end{aligned}$$

Hence the lemma holds. □

Similarly, we have

**Corollary 3.11.** *For any $z' \in \beta'[y', w'_0]$,*

$$|y'-z'| \leq a_3\, d(z').$$

*where $a_3$ is the same as in Lemma 3.10.*

Suppose that $D$ is an $a$-John domain. Then there exists an $a$-cone arc $\alpha$ in $D$ joining $x$ and $y$. Let $s_0$ bisect $\alpha$. Then

**Lemma 3.12.** *Let $u \in \alpha[x, s_0]$ and $v \in \alpha[s_0, y]$. Then for any $z \in \alpha[u, s_0]$, $d(z) \geq \frac{2\ell(\alpha[u,z])+d(u)}{4a}$, and for any $z \in \alpha[s_0, v]$, $d(z) \geq \frac{2\ell(\alpha[v,z])+d(v)}{4a}$.*

**Proof.** It suffices to prove the first statement since the proof for the second one is similar. For any $z \in \alpha[u, s_0]$, $d(z) \geq \frac{\ell(\alpha[u,z])}{a}$. If $\alpha[u, z] \subset \mathbb{B}(z, \frac{d(u)}{2})$, then $d(z) \geq \frac{d(u)}{2}$, where $\mathbb{B}(z, \frac{d(u)}{2})$ denotes the ball in $\mathbb{R}^n$ with center $z$ and radius $\frac{d(u)}{2}$. Otherwise, $d(z) \geq \frac{d(u)}{2a}$. Hence $d(z) \geq \frac{2\ell(\alpha[u,z])+d(u)}{4a}$. □

**Lemma 3.13.** (3.1) *holds for any $s_1$, $s_2 \in \alpha[x, s_0]$ (or $\alpha[x, s_0]$).*

**Proof.** It suffices to prove the first case since the proof for the second one is similar. Lemma 3.12 yields that for any $s_1, s_2 \in \alpha[x, s_0]$,

$$\begin{aligned} k_D(s_1, s_2) &\leq \ell_k(\alpha[s_1, s_2])\\ &= \int_{\alpha[s_1,s_2]} \frac{|dz|}{d(z)}\\ &\leq 4a^2 \log\Big(1+\frac{d(s_2)}{d(s_1)}\Big)\\ &\leq 4a^2 k_D(s_1, s_2)+4a^2, \end{aligned}$$

from which the proof follows. □

Let $d(v'_1) = \max\{d(u') : u' \in \alpha'[x', s'_0]\}$ and $d(v'_2) = \max\{d(u') : u' \in \alpha'[y', s'_0]\}$. Hence it follows from Lemma 3.10 and Corollary 3.11 that

**Lemma 3.14.** *For any $z' \in \alpha'[x', v'_1]$, $|x'-z'| \leq a_3\, d(z')$ and for any $z \in \alpha'[v'_1, s'_0]$, $|s'_0 - z'| \leq a_3\, d(z')$.*

Similarly,

**Corollary 3.15.** *For any* $z' \in \alpha'[y', v_2']$, $|y'-z'| \le a_3\, d(z')$ *and for any* $z \in \alpha'[v_2', s_0']$, $|s_0' - z'| \le a_3\, d(z')$.

3.3. **The proof of Theorem 1.4.** Let $z_1,\ z_2 \in D$ and $\gamma$ be a quasihyperbolic geodesic joining $z_1$, $z_2$ in $D$. In the following, we prove that $\gamma$ is a $b_1$-cone arc, that is, for any $y \in \gamma$,

$$\min\{\ell(\gamma[z_1, y]),\ \ell(\gamma[z_2, y])\} \le b_1\, d(y),$$

where $b_1 = 4a_4 e^{a_4}$, $a_4 = a_5^{2c'M}$, $a_5 = a_6^{8a^2M+C}$ and $a_6 = (8a_3)^{8c'M} e^{2C}$. It is no loss of generality to assume that $d(z_1) \le d(z_2)$.

Let $x_0 \in \gamma[z_1, z_2]$ be such that

$$d(x_0) = \max_{z\in\gamma[z_1,z_2]} d(z).$$

Then there exists an integer $t_1 \ge 0$ such that

$$2^{t_1}\, d(z_1) \le d(x_0) < 2^{t_1+1}\, d(z_1).$$

Let $y_0$ be the first point in $\gamma[z_1, x_0]$ from $z_1$ to $x_0$ with

$$d(y_0) = 2^{t_1}\, d(z_1).$$

Observe that if $d(x_0) = d(z_1)$, then $y_0 = z_1 = x_0$.

Let $y_1 = z_1$. If $z_1 = y_0$, we let $y_2 = x_0$. It is possible that $y_2 = y_1$. If $z_1 \ne y_0$, then we let $y_2, \dots, y_{t_1+1}$ be the points such that for each $i \in \{2, \dots, t_1 + 1\}$, $y_i$ denotes the first point in $\gamma[z_1, x_0]$ from $y_1$ to $x_0$ satisfying

$$d(y_i) = 2^{i-1}\, d(y_1).$$

Then $y_{t_1+1} = y_0$. We let $y_{t_1+2} = x_0$. It is possible that $y_{t_1+2} = y_{t_1+1} = x_0 = y_0$. This possibility occurs once $x_0 = y_0$.

For any fixed $i \in \{1, \dots, t_1 + 1\}$, let $\alpha_i$ be an $a$-cone arc joining $y_i$ and $y_{i+1}$ in $D$ and let $v_i$ bisect $\alpha_i$. Without loss of generality, we may assume that $d(y_i') \le d(y_{i+1}')$.

For any $z \in \alpha_i[y_i, v_i]$, Lemma 3.12 implies that

$$k_D(y_i, z) \le \ell_k(\alpha_i[y_i, z]) \le 2a \log\Big(1 + \frac{2\ell(\alpha_i[y_i, z])}{d(y_i)}\Big). \tag{3.16}$$

Similarly, for any $z \in \alpha_i[y_{i+1}, v_i]$, we have

$$k_D(y_{i+1}, z) \le 2a \log\Big(1 + \frac{2\ell(\alpha_i[y_{i+1}, z])}{d(y_{i+1})}\Big).$$

Hence

$$
\begin{aligned}
(3.17)\qquad k_D(y_i,y_{i+1}) &\leq k_D(y_i,v_i)+k_D(y_{i+1},v_i)\\
&\leq 2a\Big(\log\Big(1+\frac{2\ell(\alpha_i[y_{i+1},v_i])}{d(y_{i+1})}\Big)\\
&\quad +\log\Big(1+\frac{2\ell(\alpha_i[y_i,v_i])}{d(y_i)}\Big)\Big)\\
&\leq 4a\log\Big(1+\frac{\ell(\alpha_i)}{d(y_i)}\Big).
\end{aligned}
$$

**Lemma 3.18.** $k_D(y_i,y_{i+1})\leq a_4$.

**Proof.** Suppose that

$$
(3.19)\qquad k_D(y_i,y_{i+1})>a_4.
$$

Then

$$
\begin{aligned}
a_4 &< k_D(y_i,y_{i+1})\\
&\leq Mk_{D'}(y'_i,y'_{i+1})+C\\
&\leq c'M\log\Big(1+\frac{|y'_i-y'_{i+1}|}{d(y'_i)}\Big)+C,
\end{aligned}
$$

which implies that

$$
(3.20)\qquad |y'_i-y'_{i+1}|\geq a_5 d(y'_i).
$$

**Claim 3.21.** $d(y_i)<\lambda_D(y_i,y_{i+1})$.

Otherwise, $[y_i,y_{i+1}]\subset\mathbb{B}(y_i,d(y_i))\cap\mathbb{B}(y_{i+1},d(y_{i+1}))$, which implies that $k_D(y_i,y_{i+1})<2$. This contradicts with (3.19). Hence Claim 3.21 holds.

**Claim 3.22.** $\ell(\alpha_i)\geq a_5\lambda_D(y_i,y_{i+1})$.

Suppose not. Then (3.17) yields

$$
\begin{aligned}
(3.23)\qquad \frac{\ell(\gamma[y_i,y_{i+1}])}{2d(y_i)} &\leq \ell_k(\gamma[y_i,y_{i+1}])\\
&= k_D(y_i,y_{i+1})\\
&\leq 4a\log\Big(1+\frac{\ell(\alpha_i)}{d(y_i)}\Big)\\
&\leq 4a\log\Big(1+\frac{a_5\lambda_D(y_i,y_{i+1})}{d(y_i)}\Big).
\end{aligned}
$$

A necessary condition for (3.23) is

$$
(3.24)\qquad \lambda_D(y_i,y_{i+1})\leq a_5^2\, d(y_i).
$$

Hence (3.23) implies that $k_D(y_i, y_{i+1}) \leq a_4$. This contradiction shows that Claim 3.22 holds.

By Claim 3.22, we have

$$d(v_i) \geq \frac{\ell(\alpha_i)}{2a} \geq \frac{a_5}{2a}\lambda_D(y_i, y_{i+1}) > a_6\, \lambda_D(y_i, y_{i+1}).$$

Then Claim 3.21 guarantees that there exists $v_{i,0} \in \alpha_i[y_i, v_i]$ such that

$$d(v_{i,0}) = a_6\, \lambda_D(y_i, y_{i+1}). \tag{3.25}$$

Since by (3.16),

$$\begin{aligned} k_D(y_i, v_{i,0}) &\leq 4a \log\Big(1 + \frac{\ell(\alpha_i[y_i, v_{i,0}])}{d(y_i)}\Big) \\ &\leq 4a \log\Big(1 + \frac{ad(v_{i,0})}{d(y_i)}\Big) \\ &\leq 4a^2 a_6 \log\Big(1 + \frac{\lambda_D(y_i, y_{i+1})}{d(y_i)}\Big), \end{aligned}$$

we infer from (3.19) and the similar reasoning as in the proof of Claim 3.22 that

$$k_D(y_i, v_{i,0}) \leq \frac{1}{a_5} k_D(y_i, y_{i+1}).$$

By Claim 3.21 and (3.25),

$$k_D(y_i, v_{i,0}) \geq \log \frac{d(v_{i,0})}{d(y_i)} \geq \log a_6 > C.$$

Thus (3.16) and (3.20) imply that

$$\begin{aligned} \log\Big(1 + \frac{|y_i' - v_{i,0}'|}{d(y_i')}\Big) &\leq k_{D'}(y_i', v_{i,0}') \\ &\leq M k_D(y_i, v_{i,0}) + C \\ &< 2M k_D(y_i, v_{i,0}) \\ &< \frac{2M}{a_5} k_D(y_i, y_{i+1}) \\ &\leq \frac{2M^2}{a_5} k_{D'}(y_i', y_{i+1}') + \frac{2CM}{a_5} \\ &\leq \frac{4c'M^2}{a_5} \log\Big(1 + \frac{|y_i' - y_{i+1}'|}{d(y_i')}\Big) \\ &< \log\Big(1 + \frac{|y_i' - y_{i+1}'|}{a_5 d(y_i')}\Big). \end{aligned}$$

Hence

$$(3.26)\qquad |y_i' - v_{i,0}'| < \frac{1}{a_5}|y_i' - y_{i+1}'|,$$

which together with (3.20) give

$$(3.27)\qquad d(v_{i,0}') \le |y_i' - v_{i,0}'| + d(y_i') \le \frac{2}{a_5}|y_i' - y_{i+1}'|.$$

**Claim 3.28.** $|y_i' - v_i'| < \frac{|y_i'-y_{i+1}'|}{2}$.

Suppose for the contrary that $|y_i' - v_i'| \ge \frac{|y_i'-y_{i+1}'|}{2}$. Let $u_{0,i}' \in \gamma'[y_i', y_{i+1}']$ satisfying $d(u_{0,i}') = \max\{d(w') : w' \in \gamma'[y_i', y_{i+1}']\}$. Obviously, $\max\{|y_{i+1}' - u_{0,i}'|, |u_{0,i}' - y_i'|\} \ge \frac{|y_i'-y_{i+1}'|}{2}$. Then we know from Lemma 3.10 and Corollary 3.11 that

$$(3.29)\qquad d(u_{0,i}') \ge \frac{|y_i' - y_{i+1}'|}{2a_3}.$$

Hence by Lemma 3.10 and (3.20), there must exist some point $y_{0,i}' \in \gamma'[y_i', u_{0,i}']$ satisfying

$$(3.30)\qquad d(y_{0,i}') = \frac{|y_i' - y_{i+1}'|}{2a_3} \text{ and } |y_i' - y_{0,i}'| \le a_3\, d(y_{0,i}').$$

Since $|v_i' - v_{i,0}'| \ge |v_i' - y_i'| - |v_{i,0}' - y_i'|$, Lemma 3.14, (3.20) and (3.26) show that if $v_0' \in \alpha_i'[y_i', v_i']$ satisfies $d(v_0') = \max\{d(u') : u' \in \alpha_i'[y_i', v_i']\}$ then it must lie in $\alpha_i'[v_{i,0}', v_i']$.

Obviously, $\max\{|v_i' - v_0'|, |v_0' - y_i'|\} \ge \frac{|y_i'-y_{i+1}'|}{4}$. We know from Lemma 3.14 and Corollary 3.15 that $d(v_0') \ge \frac{|y_i'-y_{i+1}'|}{4a_3}$. By (3.27) and Lemma 3.14, we see that there exists some point $u_0' \in \alpha_i'[v_{i,0}', v_i']$ such that

$$(3.31)\qquad d(u_0') = \frac{|y_i' - y_{i+1}'|}{4a_3} \text{ and } |y_i' - u_0'| \le a_3\, d(u_0').$$

Hence (3.30) shows that

$$\begin{aligned}
\log\frac{d(u_0)}{d(y_{0,i})} &\le k_D(y_{0,i}, u_0)\\
&\le Mk_{D'}(y_{0,i}', u_0') + C\\
&\le Mc'\log\Big(1 + \frac{|u_0' - y_{0,i}'|}{\min\{d(u_0'), d(y_{0,i}')\}}\Big) + C\\
&\le Mc'\log\Big(1 + \frac{|u_0' - y_i'| + |y_i' - y_{0,i}'|}{\min\{d(u_0'), d(y_{0,i}')\}}\Big) + C\\
&\le Mc'\log\Big(1 + \frac{a_3 d(u_0') + |y_i' - y_{0,i}'|}{\min\{d(u_0'), d(y_{0,i}')\}}\Big) + C\\
&< Mc'\log(1 + 3a_3) + C,
\end{aligned}$$

which yields that

$$d(u_0) \leq (1+3a_3)^{Mc'} e^C d(y_{0,i}). \tag{3.32}$$

Lemma 3.12, (3.27) and (3.31) make sure that

$$\begin{aligned} 4a^2 M \log\Big(1+\frac{d(u_0)}{d(v_{i,0})}\Big)+C &\geq M\ell_k(\alpha_i[v_{i,0},u_0])+C \\ &\geq Mk_D(v_{i,0},u_0)+C \\ &\geq k_{D'}(v'_{i,0},u'_0) \\ &\geq \log\frac{d(u'_0)}{d(v'_{i,0})} \\ &\geq \log\frac{a_5}{8a_3}, \end{aligned}$$

whence $d(u_0) \geq a_6 d(v_{i,0})$. So we infer from Claim 3.21 and (3.25) that

$$d(u_0) \geq a_6 d(v_{i,0}) = a_6^2 \lambda_D(y_i, y_{i+1}) \geq \frac{a_6^2}{2} d(y_{i+1}) \geq \frac{a_6^2}{2} d(y_{0,i}),$$

which contradicts with (3.32). Hence Claim 3.28 holds.

It is obvious from Claim 3.28 that $|y'_{i+1}-v'_i| > \frac{|y'_i-y'_{i+1}|}{2}$. Let $q'_0 \in \alpha'_i[y'_i,v'_i] \cap \mathbb{S}(v'_i, \frac{|y'_i-v'_i|}{4a_3})$ and $u'_1 \in \alpha'_i[y'_{i+1},v'_i] \cap \mathbb{S}(v'_i, \frac{|y'_i-v'_i|}{4a_3})$. By Lemma 3.14 and Corollary 3.15, we get

$$d(q'_0) \geq \frac{|y'_i-v'_i|}{4a_3^2} \text{ and } d(u'_1) \geq \frac{|y'_i-v'_i|}{4a_3^2}. \tag{3.33}$$

Hence we have

$$\begin{aligned} \Big|\log\frac{d(u_1)}{d(q_0)}\Big| &\leq k_D(u_1,q_0) \\ &\leq Mk_{D'}(u'_1,q'_0)+C \\ &\leq Mc'\log\Big(1+\frac{|u'_1-q'_0|}{\min\{d(q'_0),d(u'_1)\}}\Big)+C \\ &\leq Mc'\log\Big(1+\frac{|u'_1-v'_i|+|v'_i-q'_0|}{\min\{d(q'_0),d(u'_1)\}}\Big)+C \\ &\leq Mc'\log(1+2a_3)+C, \end{aligned} \tag{3.34}$$

which implies that

$$\frac{d(u_1)}{(1+2a_3)^{Mc'}e^C} \leq d(q_0) \leq (1+2a_3)^{Mc'} e^C d(u_1). \tag{3.35}$$

**Claim 3.36.** $d(q_0) \geq a_5 d(v_{i,0})$.

Otherwise, Lemma 3.12, (3.25), (3.34) and (3.35) show that

$$\begin{aligned}
(3.37)\qquad \frac{\ell(\gamma[y_i,y_{i+1}])}{2d(y_i)} &\leq \ell_k(\gamma[y_i,y_{i+1}]) \\
&= k_D(y_i,y_{i+1}) \\
&\leq k_D(y_i,q_0)+k_D(q_0,u_1)+k_D(u_1,y_{i+1}) \\
&\leq 4a^2\log\Big(1+\frac{d(q_0)}{d(y_i)}\Big)+Mc'\log\Big(1+2a_3\Big) \\
&\quad +C+4a^2\log\Big(1+\frac{d(u_1)}{d(y_{i+1})}\Big) \\
&\leq 9a^2a_5\log\Big(1+\frac{\lambda_D(y_i,y_{i+1})}{d(y_i)}\Big).
\end{aligned}$$

A necessary condition for (3.37) is $\lambda_D(y_i,y_{i+1})\leq a_5^2 d(y_i)$. Hence by (3.37), we know that

$$k_D(y_i,y_{i+1})\leq 9a^2a_5\log(1+a_5^2),$$

which contradicts with (3.19). We complete the proof of Claim 3.36.

By (3.20) and (3.29)

$$|u'_{0,i}-y'_i|\geq d(u'_{0,i})-d(y'_i)\geq \frac{1}{3a_3}|y'_{i+1}-y'_i|.$$

Then Claim 3.28 guarantees that there exists $y'_0\in\gamma'[y'_i,u'_{0,i}]$ such that

$$\frac{|y'_i-v'_i|}{2a_3}=|y'_0-y'_i|.$$

Hence Lemma 3.10 implies that

$$\frac{|y'_i-v'_i|}{2a_3}=|y'_0-y'_i|\leq a_3 d(y'_0).$$

Hence (3.33) gives

$$\begin{aligned}
\log\frac{d(q_0)}{d(y_0)} &\leq k_D(q_0,y_0) \\
&\leq Mk_{D'}(q'_0,y'_0)+C \\
&\leq Mc'\log\Big(1+\frac{|y'_0-q'_0|}{\min\{d(q'_0),d(y'_0)\}}\Big)+C \\
&\leq Mc'\log\Big(1+\frac{|y'_i-v'_i|+|v'_i-q'_0|+|y'_i-y'_0|}{\min\{d(q'_0),d(y'_0)\}}\Big)+C \\
&\leq Mc'\log(1+3a_3+4a_3^2)+C.
\end{aligned}$$

We infer from Claim 3.21 and (3.25) that

$$\begin{aligned} d(q_0) &\le (1+3a_3+4a_3^2)^{Mc'}e^C d(y_0) \\ &\le 2(1+3a_3+4a_3^2)^{Mc'}e^C d(y_i) \\ &\le 2(1+3a_3+4a_3^2)^{Mc'}e^C \lambda_D(y_i, y_{i+1}) \\ &= \frac{2(1+3a_3+4a_3^2)^{Mc'}e^C}{a_6} d(v_{i,0}), \end{aligned}$$

which contradicts with Claim 3.36. We complete the proof of Lemma 3.18. □

For $i \in \{1, \cdots, t_1+1\}$, Lemma 3.18 shows

$$\frac{\ell(\gamma[y_i, y_{i+1}])}{2d(y_i)} \le \ell_k(\gamma[y_i, y_{i+1}]) = k_D(y_i, y_{i+1}) \le a_4, \tag{3.38}$$

which implies that

$$\ell(\gamma[y_i, y_{i+1}]) \le 2a_4\, d(y_i). \tag{3.39}$$

Further, for any $y \in \gamma[y_i, y_{i+1}]$, it follows from (3.38) that

$$\log \frac{d(y_i)}{d(y)} \le k_D(y, y_i) \le \; k_D(y_i, y_{i+1}) \le a_4, \tag{3.40}$$

whence

$$d(y_i) \le e^{a_4} d(y).$$

For any $y \in \gamma[y_1, x_0]$, there is some $i \in \{1, \cdots, t_1+1\}$ such that $y \in \gamma[y_i, y_{i+1}]$. It follows from (3.39) and (3.40) yield that

$$\begin{aligned} \ell(\gamma[z_1, y]) &= \ell(\gamma[y_1, y_2]) + \ell(\gamma[y_2, y_3]) + \cdots + \ell(\gamma[y_i, y]) \\ &\le 2a_4(d(y_1) + d(y_2) + \cdots + d(y_i)) \\ &\le 4a_4\, d(y_i) \\ &\le 4a_4 e^{a_4}\, d(y). \end{aligned} \tag{3.41}$$

By replacing $\gamma[z_1, x_0]$ by $\gamma[z_2, x_0]$ and repeating the procedure as above, we also get that

$$\ell(\gamma[z_2, y]) \le 4a_4 e^{a_4}\, d(y). \tag{3.42}$$

The combination of (3.41) and (3.42) conclude the proof of Theorem 1.4.

## 4. The proof of Theorem 1.6

Before the proof of Theorem 1.6, we introduce the following lemma which is the main result in [7].

**Lemma E.** ([7, Theorem 1.1]) *Suppose that $D \subset \mathbb{R}^n$ is a domain which is homeomorphic to a c-uniform domain via a quasiconformal mapping $f$. If $\gamma$ is a quasihyperbolic geodesic in $D$ and if $L$ is any other arc in $D$ with the same end points as $\gamma$, then*

$$\ell(\gamma) \leq b\, \ell(L),$$

*where $\ell(\gamma)$ denotes the length of $\gamma$ and the positive constant $b$ depends only on $n$, $c$, and the dilatation $K(f)$ of $f$.*

4.1. **The proof of Theorem 1.6.** For any $z_1$ and $z_2 \in D$, we let $\gamma$ be a quasihyperbolic geodesic joining $z_1$ and $z_2$ in $D$. It follows from Theorem 1.4 and Lemma E that

(1) $\min\limits_{j=1,2} \ell(\gamma[z_j, z]) \leq b'\, d(z)$ for all $z \in \gamma$, and

(2) $\ell(\gamma) \leq b'\, \lambda_D(z_1, z_2)$,

where the constant $b'$ depends on only $a$, $c$, $K$ and $n$. Hence the proof of Theorem 1.6 easily follows from Theorem D.

M. Huang, Department of Mathematics, Hunan Normal University, Changsha, Hunan 410081, People's Republic of China
*E-mail address*: mzhuang79@yahoo.com.cn

X. Wang, Department of Mathematics, Hunan Normal University, Changsha, Hunan 410081, People's Republic of China
*E-mail address*: xtwang@hunnu.edu.cn